\documentclass[12pt,leqno]{article}
\title{On Kummer and Stickelberger relations }
\author{Roland Qu\^eme}
\usepackage{amsmath,amsthm,amstext,amsbsy}
\newtheorem{thm}{Theorem}[section]
\newtheorem{cor}[thm]{Corollary}

\newtheorem{lem}[thm]{Lemma}
\font\mathbb=msbm10

\newcommand{\Q}{\mbox{\mathbb Q}}
\newcommand{\Z}{\mbox{\mathbb Z}}

\newcommand{\modu}{\ \mbox{mod}\ }
\newcommand{\be}{\begin{equation}}
\newcommand{\ee}{\end{equation}}
\date{2005 dec 30}
\begin{document}
\maketitle
\abstract
Roland Qu\^eme

13 avenue du ch\^ateau d'eau

31490 Brax

France

2005 dec 30

tel : 0561067020

cel : 0684728729

mailto: roland.queme@wanadoo.fr

home page: http://roland.queme.free.fr/

************************************

Let $p$ be an odd  prime. Let ${\bf F}_p$ be the finite field of $p$ elements with no null part ${\bf F}_p^*$. Let $K_p=\Q(\zeta_p)$ be the $p$-cyclotomic field. 
Let $\pi$ be the prime ideal of $K_p$ lying over $p$.
Let $v$ be a primitive root $\modu p$.
In the sequel of this paper, for $n\in \Z$  let us note briefly  $v^n$ for $v^n \modu p$ with $1\leq v^n\leq p-1$. For $n<0$ it is understood with $v^n\times v^{-n}\equiv 1\modu p$.  
Let $\sigma :\zeta_p\rightarrow \zeta_p^v$ be a $\Q$-isomorphism of $K_p/\Q$.  
Let $G_p$ be the Galois group of $K_p/\Q$.
Let $P(\sigma)=\sum_{i=0}^{p-2}\sigma^i\times v^{-i},\quad P(\sigma)\in\Z[G_p]$.
Let $A\in\Z[\zeta_p]$  with $\pi^{2m+1} \ \|\  A-1$ where   $m$ is a natural number $1\leq m\leq \frac{p-3}{2}$. 
It is possible to show in an elementary way that $A^{P(\sigma)}\equiv 1\modu \pi^{p-1}$.

We suppose that $p$ is an irregular prime and that $p$ does not divide the class number of the maximal totally real subfield $K_p^+$ of $K_p$. 
Let $C_p$ be the  $p$-class group of $K_p$. Let $\Gamma$ be a subgroup of $C_p$ of order $p$ annihilated by $\sigma-\mu$ with 
$\mu\in{\bf F}_p^*$.  
From Kummer, there exists    not principal prime ideals  $\mathbf q$  of $\Z[\zeta_p]$ of inertial degree $1$ with class $Cl(\mathbf q)\in \Gamma$.
There exist     singular {\bf not}  primary numbers $A$ with $A \Z[\zeta_p]= \mathbf q^p$ and $\pi^{2m+1}\ \|\ A-1$  where $m$ and $\mu$ are connected by 
$\mu\equiv v^{2m+1}\modu p,\ 1\leq m\leq \frac{p-3}{2}$.
We prove,  by an application of Stickelberger's relation to the prime ideal $\mathbf q$,   that now  we can {\it climb} up to the $\pi$-adic congruence 
$A^{P(\sigma)}\equiv 1\modu \pi^{2(p-1)}$.

This $\pi$-adic improvement allows us to {\it caracterize} in a straightforward way the $p$-class group $C_p$ by the congruence $\modu p$. 
With $v,m$ defined above:
\begin{equation} 
\sum_{i=1}^{p-2} v^{(2m+1)(i-1)}\times(\frac{v^{-(i-1)}-v^{-i}\times v}{p})\equiv 0\modu p:
\end{equation}
The numerical verification of this  congruence is completely consistent with table of irregular primes in Washington \cite{was} p. 410.

%

%
\section{Some definitions}
\begin{enumerate}
\item
Let $p$ be an odd prime. Let $\zeta_p$ be a root of the polynomial equation $X^{p-1}+X^{p-2}+\dots+X+1=0$.
Let $K_p$ be the $p$-cyclotomic field $K_p=\Q(\zeta_p)$. The ring of integers of $K_p$ is $\Z[\zeta_p]$.
Let $K_p^+$ be the maximal totally real subfield of $K_p$.
Let $v$ be a primitive root $\modu p$ and $\sigma: \zeta_p\rightarrow \zeta_p^v$ be a $\Q$-isomorphism of $K_p$.
Let $G_p$ be the Galois group of $K_p/\Q$.
Let ${\bf F}_p$ be the finite field of cardinal $p$ with no null part  ${\bf F}_p^*$.
Let $\lambda=\zeta_p-1$. The prime ideal of $K_p$ lying over $p$ is $\pi=\lambda \Z[\zeta_p]$.
\item 
Suppose that $p$ is irregular and that $p$ does not divide  the class number of $K_p^+$. 
Let $C_p$ be the $p$-class group of $K_p$.
Let $r$ be the rank of $C_p$.
Let $\Gamma$ be a  subgroup of order $p$ of $C_p$ annihilated by $\sigma-\mu\in {\bf F}_p[G_p]$ with $\mu\in{\bf F}_p^*$. Then $\mu\equiv v^{2m+1}\modu p$ with a natural integer $m,\quad 1\leq m\leq \frac{p-3}{2}$.
\item 
An  integer  $A\in \Z[\zeta_p]$ is said singular if $A^{1/p}\not\in K_p$ and if  there exists an   ideal $\mathbf a$  of $\Z[\zeta_p]$ such that 
$A \Z[\zeta_p]=\mathbf a^p$.
There exists  singular integers $A$  with $A \Z[\zeta_p] =\mathbf a^p$ where $\mathbf a$ is a {\bf not} principal  ideal of $\Z[\zeta_p]$ 
verifying   
\begin{equation}\label{e512101}
\begin{split}
& Cl(\mathbf a)\in \Gamma,\\
& \sigma(A)=A^\mu\times\alpha^p,\\
\end{split}
\end{equation}
with $\alpha\in K_p$ and $\pi^{2m+1} \ \|\ A-1$.
Moreover, the number $A$ verifies \begin{equation}\label{e512103}
A\times\overline{A}=D^p,
\end{equation}
for some integer $D\in O_{K_p^+}$. 
\end{enumerate} 
%
\section{On Kummer and Stickelberger relation}
\begin{enumerate}
\item
Here we fix a notation for the sequel. 
Let $v$ be a primitive root $\modu p$. For every integer $k\in\Z$   then  $v^k$ is understood $\modu p$ so $1\leq v^k\leq p-1$.
 If $k<0$ it is to be understood as $v^kv^{-k}\equiv 1\modu p$.

\item
Let $q\not=p$ be an odd prime. 
Let $\zeta_q$ be a root of the minimal polynomial equation $X^{q-1}+X^{q-2}+\dots+X+1=0$.
Let $K_q=\Q(\zeta_q)$ be the $q$-cyclotomic field.
The ring of integers of $K_q$ is $\Z[\zeta_q]$.
Here we fix a notation for the sequel. 
Let $u$ be a primitive root $\modu q$. For every integer $k\in\Z$   then  $u^k$ is understood $\modu q$ so $1\leq u^k\leq q-1$.
 If $k<0$ it is to be understood as $u^ku^{-k}\equiv 1\modu q$.
Let $K_{pq}=\Q(\zeta_p,\zeta_q)$. Then $K_{pq}$ is the compositum $K_pK_q$.
The ring of integers of $K_{pq}$ is $\Z[\zeta_{pq}]$.
\item
Let $\mathbf q$ be a prime ideal of $\Z[\zeta_p]$ lying over the prime $q$.
Let $m=N_{K_p/\Q}(\mathbf q)= q^f$ where $f$ is the smallest integer such that $q^f\equiv 1\modu p$.
If $\psi(\alpha)=a$ is the image of $\alpha\in \Z[\zeta_p]$ under the natural map 
$\psi: \Z[\zeta_p]\rightarrow \Z[\zeta_p]/\mathbf q$, then for 
$\psi(\alpha)=a\not\equiv 0$ define a character $\chi_{\mathbf q}^{(p)}$ on ${\bf F}_m=\Z[\zeta_p]/\mathbf q$ by
\begin{equation}
\chi_{\mathbf q}^{(p)}(a)=\{\frac{\alpha}{\mathbf q}\}_p^{-1}=\overline{\{\frac{\alpha}{\mathbf q}\}}_p,
\end{equation}
where $\{\frac{\alpha}{\mathbf q}\}=\zeta_p^c$ for some natural integer $c$,  is the $p^{th}$ power residue character $\modu \mathbf q$.
We define 
\begin{displaymath}
g(\mathbf q)=\sum_{x\in{\bf F}_m}(\chi_{\mathbf q}^{(p)}(x)\times 
\zeta_q^{Tr_{{\bf F}_m/{\bf F}_q}(x)})\in \Z[\zeta_{pq}],
\end{displaymath}
and $\mathbf G(\mathbf q)= g(\mathbf q)^p$.
It follows that $G(\mathbf q)\in \Z[\zeta_{pq}]$.
Moreover $\mathbf G(\mathbf q)=g(\mathbf q)^p\in \Z[\zeta_p]$, see for instance Mollin \cite{mol} prop. 5.88 (c) p. 308
or Ireland-Rosen \cite{ire} prop. 14.3.1 (c) p. 208.
\end{enumerate}
%
The Stickelberger's relation is classically:
\begin{thm}\label{t12201}

In $\Z[\zeta_p]$ we have the ideal decomposition 
\begin{equation}\label{e512121}
\mathbf G(\mathbf q)\Z[\zeta_p]=\mathbf q^{S},
\end{equation}
with $S=\sum_{t=1}^{p-1} t\times \varpi_t^{-1}$
where  $\varpi_t\in Gal(K_p/\Q)$ is given by $\varpi_t: \zeta_p\rightarrow \zeta_p^t$.
\end{thm}
See for instance Mollin \cite{mol} thm. 5.109 p. 315 and Ireland-Rosen \cite{ire} thm. 2. p.209.
%
\subsection{On the structure of $\mathbf G(\mathbf q)$.}
In this subsection we are studying carefully the structure of $g(\mathbf q)$ and $\mathbf G(\mathbf q)$.
\begin{lem}\label{l512151}
If $q\not\equiv 1\modu p$ then $g(\mathbf q)\in \Z[\zeta_p]$.
\begin{proof}$ $
\begin{enumerate}
\item
Let $u$ be a primitive root $\modu q$. Let $\tau :\zeta_q\rightarrow \zeta_q^u$ be a $\Q$-isomorphism generating $Gal(K_q/\Q)$.
The isomorphism $\tau$ is extended to a $K_p$-isomorphism of $K_{pq}$ by 
$\tau:\zeta_q\rightarrow \zeta_q^u,\quad \zeta_p\rightarrow \zeta_p$.
Then  $g(\mathbf q)^p=\mathbf G(\mathbf q)\in \Z[\zeta_p]$ and so 
\begin{displaymath}
\tau(g(\mathbf q))^p=g(\mathbf q)^p,
\end{displaymath}
and it follows that there exists a natural integer $\rho$ with $\rho<p$ such that 
\begin{displaymath}
\tau(g(\mathbf q))= \zeta_p^\rho\times  g(\mathbf q).
\end{displaymath}
Then $N_{K_{pq}/K_p}(\tau(g(\mathbf q)))=\zeta_p^{(q-1)\rho}\times N_{K_{pq}/K_p}(g(\mathbf q))$ and so  $\zeta_p^{\rho(q-1)}=1$.
\item
If $q\not\equiv 1\modu p$, it implies that $\zeta_p^\rho=1$ and so that $\tau(g(\mathbf q))=g(\mathbf q)$ and thus that $g(\mathbf q)\in \Z[\zeta_p]$.
\end{enumerate}
\end{proof}
\end{lem}
%
Let us note in the sequel $g(\mathbf q)=\sum_{i=0}^{q-2} g_i\times \zeta_q^i$ with $g_i\in \Z[\zeta_p]$.
\begin{lem}\label{l512152}
If $q\equiv 1\modu p$ then $g_0=0$.
\begin{proof}
Suppose that $g_0\not=0$ and search for a contradiction:
we start of 
\begin{displaymath}
\tau(g(\mathbf q))= \zeta_p^\rho\times  g(\mathbf q).
\end{displaymath}
We have $g(\mathbf q)=\sum_{i=0}^{q-2} g_i\times \zeta_q^i$  and so 
$\tau(g(\mathbf q))=\sum_{i=0}^{q-2}  g_i\times \zeta_q^{i u}$, 
therefore 
\begin{displaymath}
\sum_{i=0}^{q-2} (\zeta_p^\rho g_i)\times \zeta_q^i=\sum_{i=0}^{q-2}  g_i\times \zeta_q^{i u},
\end{displaymath}
thus $g_0=\zeta_p^\rho \times g_0$ and so $\zeta_p^\rho=1$ which 
implies that $\tau(g(\mathbf q))=g(\mathbf q)$ and so $g(\mathbf q)\in \Z[\zeta_p]$.
Then $\mathbf G(\mathbf q)=g(\mathbf q)^p$ and so Stickelberger relation leads to 
$g(\mathbf q)^p \Z[\zeta_p] =\mathbf q^S$ where $S=\sum_{t=1}^{p-1}t\times\varpi_t^{-1}$.
Therefore $\varpi_1^{-1}(\mathbf q) \ \|\ \mathbf q^S$ because $q$ splits totally in $K_p/\Q$
 and $\varpi_t^{-1}(\mathbf q)\not=\varpi_{t^\prime}^{-1}(\mathbf q)$ for $t\not=t^\prime$. This case is not possible because the first member 
$g(\mathbf q)^p$ is a $p$-power.
\end{proof}
\end{lem}
%
\begin{lem}\label{l512152}
If $q\equiv 1\modu p$ then 
\begin{equation}\label{e512151}
\begin{split}
& \mathbf G(\mathbf q) = g(\mathbf q)^p,\\
&g(\mathbf q)=\zeta_q +\zeta_p^\rho\zeta_q^{u^{-1}}+\zeta_p^{2\rho}\zeta_q^{u^{-2}}+\dots \zeta_p^{(q-2)\rho}\zeta_q^{u^{-(q-2)}},\\
& g(\mathbf q)^p \Z[\zeta_p] =\mathbf q^S.\\
\end{split}
\end{equation}
\begin{proof} $ $
\begin{enumerate}
\item
We start of $\tau(g(\mathbf q))=\zeta_p^\rho\times g(\mathbf q)$ and so 
\begin{equation}\label{e512152}
\sum_{i=1} ^{q-2}g_i \zeta_q^{ui}=\zeta_p^\rho\times\sum_{i=1}^{q-2} g_i \zeta_q^i,
\end{equation}
which implies that $g_i=g_1\zeta_p^\rho $ for  $u\times i\equiv 1\modu q$ and so $g_{u^{-1}}=g_1\zeta_p^\rho$ (where $u^{-1}$ is to be understood by
$u^{-1}\modu q$,  so $1\leq u^{-1}\leq q-1)$.
\item
Then
$\tau^2(g(\mathbf q))=\tau(\zeta_p^{\rho} g(\mathbf q))=\zeta_p^{2\rho} g(\mathbf q)$.
Then 
\begin{displaymath}
\sum_{i=1} ^{q-2}g_i \zeta_q^{u^2i}=\zeta_p^{2\rho}\times (\sum_{i=1}^{q-2} g_i \zeta_q^i),
\end{displaymath}
which implies that $g_i=g_1\zeta_p^{2\rho}$ for $u^2\times i\equiv 1\modu q$ and so $g_{u^{-2}}=g_1\zeta_p^{2\rho}$.
\item
We continue up to 
$\tau^{(q-2)\rho}(g(\mathbf q))=\tau^{q-3}(\zeta_p^\rho g(\mathbf q))=\dots=\zeta_p^{(q-2)\rho} g(\mathbf q)$.
Then 
\begin{displaymath}
\sum_{i=1} ^{q-2}g_i \zeta_q^{u^{q-2}i}=\zeta_p^{(q-2)\rho}\times(\sum_{i=1}^{q-2} g_i \zeta_q^i),
\end{displaymath}
which implies that $g_i=g_1\zeta_p^{(q-2)\rho}$ for $u^{q-2}\times i\equiv 1\modu q$ and so $g_{u^{-(q-2)}}=g_1\zeta_p^{(q-2)\rho}$.
\item
Observe that $u$ is a primitive root $\modu q$ and so $u^{-1}$ is a primitive root $\modu q$.
Then it follows that 
$g(\mathbf q) =g_1\times (\zeta_q +\zeta_p^\rho\zeta_q^{u^{-1}}+\zeta_p^{2\rho}\zeta_q^{u^{-2}}+\dots \zeta_p^{(q-2)\rho}\zeta_q^{u^{-(q-2)}})$.
Let $U=\zeta_q +\zeta_p^\rho\zeta_q^{u^{-1}}+\zeta_p^{2\rho}\zeta_q^{u^{-2}}+\dots \zeta_p^{(q-2)\rho}\zeta_q^{u^{-(q-2)}}$.
\item
We prove now that $g_1\in \Z[\zeta_p]^*$. 
From Stickelberger relation $g_1^p \times U^p =\mathbf q^{S}$.
From $S=\sum_{i=1}^{p-1}\varpi_t^{-1}\times t$ it follows that 
$\varpi_t^{-1}(\mathbf q)^t\ \|\ \mathbf q^{S}$
and so that $g_1\not\equiv 0\modu \varpi_t^{-1}(\mathbf q)$
because $g_1^p$ is a $p$-power,
which implies that 
$g_1\in \Z[\zeta_p]^*$. We can take $g_1=1$ because $\mathbf G(\mathbf q)$ is defined modulo the unit group $\Z[\zeta_p]^*$.
\item
From Stickelberger,  $g(\mathbf q)^p \Z[\zeta_p]=\mathbf q^S$,
which achieves the proof.
\end{enumerate}
\end{proof}
\end{lem}
%
\paragraph{Remark:}
From
\begin{equation}
\begin{split}
&g(\mathbf q)=\zeta_q +\zeta_p^\rho\zeta_q^{u^{-1}}+\zeta_p^{2\rho}\zeta_q^{u^{-2}}+\dots +\zeta_p^{(q-2)\rho}\zeta_q^{u^{-(q-2)}},\\
&\Rightarrow \tau(g(\mathbf q))=\zeta_q^u +\zeta_p^\rho\zeta_q+\zeta_p^{2\rho}\zeta_q^{u^{-1}}+\dots +\zeta_p^{(q-2)\rho}\zeta_q^{u^{-(q-3)}},\\
&\Rightarrow\zeta^\rho\times g(\mathbf q)=\zeta^\rho\zeta_q +\zeta_p^{2\rho}\zeta_q^{u^{-1}}+\zeta_p^{3\rho}\zeta_q^{u^{-2}}+\dots +\zeta_p^{(q-1)\rho}
\zeta_q^{u^{-(q-2)}}\\
\end{split}
\end{equation}
and we can verify directly that $\tau(g(\mathbf q))=\zeta_p^\rho \times g(\mathbf q)$ for this expression of $g(\mathbf q)$, observing that $q-1\equiv 0\modu p$.
%
%
\begin{lem}\label{l12161}
Let $S=\sum_{t=1}^{p-1} \varpi_t^{-1}\times t$ where $\varpi_t$ is the $\Q$-isomorphism 
given by $\varpi_t:\zeta_p\rightarrow \zeta_p^t$ of $K_p$.
Let $v$ be a primitive root $\modu p$. Let $\sigma$ be the $\Q$-isomorphism of $K_p$ given by  $\zeta_p\rightarrow\zeta_p^v$.
Let $P(\sigma)=\sum_{i=0}^{p-2} \sigma^i\times v^{-i}\in\Z[G_p]$.
Then $S=P(\sigma)$.
\begin{proof}
Let us consider one term $\varpi_t^{-1} \times t$. Then $v^{-1}=v^{p-2}$ is a primitive root $\modu p$ because $p-2$ and $p-1$ are coprime and so there exists one and one $i$ such that 
$t=v^{-i}$. Then $\varpi_{v^{-i}}:\zeta_p\rightarrow \zeta_p^{v^{-i}}$ and so $\varpi_{v^{-i}}^{-1}:\zeta_p\rightarrow\zeta_p^{v^i}$
and so $\varpi_{v^{-i}}^{-1}=\sigma^i$ (observe that $\sigma^{p-1}\times v^{-(p-1)}=1$), which achieves the proof.
\end{proof}
\end{lem} 
%
\paragraph{Remark} : The previous lemma is  a verification of  the consistency of results in Ribenboim \cite{rib} p. 118, of Mollin \cite {mol} p. 315 and of Ireland-Rosen p. 209 with our computation. In the sequel we use Ribenboim notation more adequate for the factorization in ${\bf F}_p[G]$.
%
In that case the Stickelberger's relation is connected with the Kummer's relation on Jacobi resolvents, see for instance 
Ribenboim, \cite{rib} (2A) b. p. 118 and (2C) relation (2.6) p. 119.

\begin{lem}\label{l512162}$ $

\begin{enumerate}
\item
$g(\mathbf q)$ defined in relation (\ref{e512151}) is the  Jacobi resolvent $<\zeta_p^\rho,\zeta_q>$.
\item
$\rho=-v$.
\end{enumerate}
\begin{proof}$ $
\begin{enumerate}
\item
Apply formula of Ribenboim \cite{rib} (2.2) p. 118 with 
$p=p, q=q, \zeta =\zeta_p,\quad \rho=\zeta_q, \quad n=\rho,\quad u=i,\quad m=1$ and $h=u^{-1}$
(where the left members notations $p, q, \zeta,\rho, n, u, m$ and $h$ are the Ribenboim notations).
\item
We start of 
$<\zeta_p^\rho,\zeta_q>=g(\mathbf q)$.
Then $v$ is a primitive root $\modu p$, so there exists  a natural integer $l$ such that 
$\rho \equiv v^l\modu p$.
By conjugation $\sigma^{-l}$ we get 
$<\zeta_p,\zeta_q>=g(\mathbf q)^{\sigma^{-l}}$.
Raising to $p$-power 
$<\zeta_p,\zeta_q>^p=g(\mathbf q)^{p\sigma^{-l}}$.
From lemma \ref{l12161} and Stickelberger relation
$<\zeta_p,\zeta_q>^p=\mathbf q^{P(\sigma)\sigma^{-l}}$.
From Kummer's relation (2.6) p. 119 in Ribenboim \cite{rib}, we get 
$<\zeta_p,\zeta_q>^p=\mathbf q^{P_1(\sigma)}$ with $P_1(\sigma)=\sum_{j=0}^{p-2}\sigma^j v^{(p-1)/2-j}$.
Therefore
$\sum_{i=0}^{p-2} \sigma^{i-l}v^{-i}=\sum_{j=0}^{p-2}\sigma^j v^{(p-1)/2-j}$.
Then 
$i-l\equiv j\modu p$ and $-i\equiv \frac{p-1}{2}-j\modu p$ (or  $i\equiv j-\frac{p-1}{2}\modu p$)
imply that 
$j-\frac{p-1}{2}-l\equiv j\modu p$,
so 
$l+\frac{p-1}{2}\equiv 0\modu p$,
so
$l\equiv -\frac{p-1}{2}\modu p$,
and 
$l\equiv \frac{p+1}{2}\modu p$,
thus 
$\rho\equiv v^{(p+1)/2}\modu p$ and finally $\rho= -v$.
\end{enumerate}
\end{proof}
\end{lem}
\paragraph{Remark}: The previous lemma allows to verify the consistency of our computation with Jacobi resultents used in Kummer (see Ribenboim 
p. 118-119).
%
\begin{lem}\label{l512161}
$g(\mathbf q)\equiv -1\modu \pi$.
\begin{proof}
From $g(\mathbf q)=\zeta_q +\zeta_p^\rho\zeta_q^{u^{-1}}+\zeta_p^{2\rho}\zeta_q^{u^{-2}}+\dots +\zeta_p^{(q-2)\rho}\zeta_q^{u^{-(q-2)}}$, we see that 
$g(\mathbf q)\equiv \zeta_q +\zeta_q^{u^{-1}}+\zeta_q^{u^{-2}}+\dots +\zeta_q^{u^{-(q-2)}}\modu \pi$.
From $u^{-1}$ primitive root $\modu p$ it follows  that 
$1+\zeta_q +\zeta_q^{u^{-1}}+\zeta_q^{u^{-2}}+\dots +\zeta_q^{u^{-(q-2)}}=0$, which leads to the result.
\end{proof}
\end{lem}
%
\subsection{Singular integers and Kummer-Stickelberger's relation}
We know that the group of ideal classes of the cyclotomic field is generated by the ideal classes of prime ideals of degree $1$, see for instance 
Ribenboim, \cite{rib} (3A) p. 119.
From now we suppose  that the prime ideal $\mathbf q$ of $\Z[\zeta_p]$ has a class $Cl(\mathbf q)\in \Gamma$ where $\Gamma$ is  
a subgroup of order $p$ of $C_p$ previously defined,  with a singular integer $A$  given  by $A \Z[\zeta_p] = \mathbf q^p$ and  
$N_{K_p/\Q}(\mathbf q)=q\equiv 1\modu p$.
%
\begin{thm}\label{l512163}$ $

$(\frac{g(\mathbf q}{\overline{g(\mathbf q})})^{p^2}=(\frac{A^2}{D^p})^{P(\sigma)}$.
\begin{proof}
We start of 
$\mathbf G(\mathbf q)\Z[\zeta_p]= g(\mathbf q)^p \Z[\zeta_p]= \mathbf q^S$.
Raising to $p$-power we get 
$g(\mathbf q)^{p^2} \Z[\zeta_p]= \mathbf q^{pS}$.
But $A \Z[\zeta_p] = \mathbf q^p$, so 
$g(\mathbf q)^{p^2} \Z[\zeta_p]= A^{S}$, so 
 $g(\mathbf q)^{p^2}\times  \zeta_p^w\times \eta= A^{S},\quad \eta\in \Z[\zeta+\zeta^{-1}]^*$ where $w$ is a natural number.
Therefore, by complex conjugation, we get 
$ \overline{g(\mathbf q)}^{p^2}\times\zeta_p^{-w}\times  \eta= \overline{A}^{S}.$
Then 
$ (\frac{g(\mathbf q)}{\overline{g(\mathbf q)}})^{p^2}\times\zeta_p^{2w}=(\frac{A}{\overline{A}})^{S}$.
From $A\times\overline{A}=D^p$ we obtain
$(\frac{g(\mathbf q}{\overline{g(\mathbf q}})^{p^2}\times\zeta_p^{2w}=(\frac{A^2}{D^p})^{S}$.
We have $\pi^{2m+1}\ \|\ A-1$  for $m\geq 1$ and $\pi^{p+1}\ |\ D^p-1$ and so $(\frac{A^{2}}{D^p})^{P(\sigma)}\equiv 1\modu \pi^{2m+1}$.
From $g(\mathbf q)\equiv -1\modu \pi$ it follows that $\overline{g(\mathbf q)}\equiv -1\modu\pi$ and so that  $\frac{g(\mathbf q)}{\overline{g(\mathbf q)}}\equiv 1\modu \pi$ and so 
$(\frac{g(\mathbf q)}{\overline{g(\mathbf q)}})^p\equiv 1\modu \pi^p$ and 
thus $w=0$,  
so $(\frac{g(\mathbf q)}{\overline{g(\mathbf q)}})^{p^2}=
(\frac{A^2}{D^p})^{P(\sigma)}$.
\end{proof}
\end{thm}
%
\begin{lem}\label{l512171}
\begin{equation}
P(\sigma)=\sum_{i=0}^{p-2}\sigma^i\times v^{-i}=v^{-(p-2)}\times \{\prod_{k=0,\ k\not=1}^{p-2}(\sigma-v^{k})\}+p\times R(\sigma),
\end{equation}
where $R(\sigma)\in\Z[G_p]$ with $deg(R(\sigma))<p-2$.
\begin{proof}
Let us consider the polynomial $R_0(\sigma)=P(\sigma)-v^{-(p-2)}\times \{\prod_{k=0,\ k\not=1}^{p-2}(\sigma-v^{k})\}$  in ${\bf F}_p[G_p]$.
Then $R_0(\sigma)$ is of degree smaller than $p-2$ and  the two 
polynomials $\sum_{i=0}^{p-2} \sigma^iv^{-i} $ and  $\prod_{k=0,\ k\not=1}^{p-2}(\sigma-v^{k})$ take a null value in ${\bf F}_p[G_p]$  when $\sigma$ takes  the $p-2$ different  values  $\sigma=v^k$ for $k=0,\dots, p-2,\quad k\not= 1$. Then $R_0(\sigma)=0$ in ${\bf F}_p[G_p]$ which leads to the result in 
$\Z[G_p]$.
\end{proof}
\end{lem}
Let us note in the sequel  
\begin{equation}\label{e12201}
T(\sigma)=v^{-(p-2)}\times\prod_{k=0,\ k\not=1}^{p-2}(\sigma-v^{k}).
\end{equation}
%
\begin{thm}\label{l512164}$ $

$A^{P(\sigma)}\equiv 1\modu\pi^{2(p-1)}$.
\begin{proof}$ $ 
\begin{enumerate}
\item
From $A\times \overline{A}=D^p$ we know that $D\equiv 1\modu\pi$.
Let us consider $D\equiv 1+d\lambda^n\modu \pi^{n+1}$ where $d\in \Z,\quad d\not\equiv 0\modu p$ and $n\geq 1$. From  $D\in \Z[\zeta_p+\zeta_p^{-1}]$, we derive that $n$ is even and $n\geq 2$.
Then $\sigma(D)\equiv 1+d\sigma(\lambda)^n\equiv 1+d(\zeta^v-1)^n\equiv 1+d((\lambda+1)^v-1)^n$
$\equiv 1+d v^n\lambda^n\modu\pi^{n+1}$.   
Also $D^{v^n}=(1+d\lambda^n)^{v^n}\equiv 1+dv^n\lambda^n\modu\pi^{n+1}$.
Then $D^{\sigma-v^n}\equiv 1\modu\pi^{n+1}$.
But $D^{\sigma-v^n}\in K_p^+$ and so $\pi^{n^{\prime}}\ \|\ D^{\sigma-v^n}-1$ with $n^\prime$ even and $n^\prime>n$.
 Using the factorization of $T(\sigma)$ given, we can climb  starting $\pi$-adically of $n^\prime$,  step by step, in the $\pi$-adic approximations to get  finally 
$D^{T(\sigma)}\equiv 1\modu\pi^{p-1}$. Then $D^{pT(\sigma)}\equiv 1\modu \pi^{2(p-1)}$. From lemma \ref{l512171}, $p P(\sigma)= pT(\sigma)+p^2 \times 
R(\sigma)$
where $R(\sigma)\in\Z[G_p]$ and so 
$D^{pP(\sigma)}\equiv 1\modu\pi^{2(p-1)}$.
\item
We have $g(\mathbf q)\equiv -1\modu \pi$ in $\Z[\zeta_{pq}]$ and so 
$g(\mathbf q)^p\equiv -1\modu\pi^{p}$ in $\Z[\zeta_{pq}]$. 
Therefore $g(\mathbf q)^p-1= \lambda^p\times  a$ with 
$a=\sum_{i=0}^{p-2}\sum_{j=0}^{q-2} a_{ij} \zeta_p^i\zeta_q^j,\quad a_{ij}\in \Z$. 
But $g(\mathbf q)^p\in \Z[\zeta_p]$ so $a\in K_p$ and so $a_{ij}=0$ when $j\not=0$.
Therefore $a=\sum_{i=0}^{p-2} a_{i0}\zeta_p^i,\quad a_{i0}\in\Z$ and so $g(\mathbf q)^p\equiv 1\modu\pi^p$ in $\Z[\zeta_p]$.
Therefore $(\frac{g(\mathbf q}{\overline{g(\mathbf q)}})^{p^2}\equiv 1\modu\pi^{2(p-1)}$ in $K_{p}$.
We have seen also that $D^{pP(\sigma)}\equiv 1\modu\pi^{2(p-1)}$ in $\Z[\zeta_p]$. Therefore from theorem \ref{l512163} $A^{2P(\sigma)}\equiv 1\modu \pi^{2(p-1)}$ and 
from  $A\equiv 1\modu \pi^{2m+1}$ it follows that $A^{P(\sigma)}\equiv 1\modu \pi^{2(p-1)}$.
\end{enumerate}
\end{proof}
\end{thm}
\paragraph{Remark:} if $C\in\Z[\zeta_p]$ is a semi-primary number  with  $C\equiv 1\modu\pi^2$
we can only assert in general that $C^{P(\sigma)}\equiv  1\modu\pi^{p-1}$. For the singular numbers $A$ considered here we assert more: 
$A^{P(\sigma)}\equiv 1\modu\pi^{2(p-1)}$.  We shall use this $\pi$-adic improvement in the sequel.
%
\clearpage
\section{Congruences caracterizing the $p$-class group}
\begin{enumerate}
\item
In this section we shall give a  polynomial $Q(\sigma)\in\Z[G_p]$ explicitly computable in a straightforward way    which {\it caracterize}
the structure of the $p$-class group of $K_p$.
\item
We know  that the    $p$-class group $C_p=\oplus_{k=1}^{r} \Gamma_k$ where $\Gamma_k$ are groups of order $p$ annihilated by $\sigma-\mu_k, 
\quad \mu_k \equiv v^{2m_k+1}\modu p,\quad 1\leq m_k\leq\frac{p-3}{2}$.
Let us consider the singular numbers $A_k,\quad k=1,\dots,r$, with $\pi^{2m_k+1} \ \|\ A_k-1$  defined  in  lemmas \ref{l512171} and \ref{l512165}.
From Kummer, the group of ideal classes of $K_p$ is generated by the classes of prime ideals of degree $1$ (see for instance Ribenboim \cite{rib} (3A) p. 119).
\end{enumerate}
%
\begin{lem}\label{l512165}
\begin{equation}\label{e512172}
P(\sigma)\times (\sigma-v)= T(\sigma)\times(\sigma-v)+pR(\sigma)\times (\sigma-v)=p\times Q(\sigma),
\end{equation}
where $Q(\sigma)=\sum_{i=1}^{p-2}\delta_i\times \sigma^i\in\Z[G_p]$ is given by  
\begin{equation}
\begin{split}
& \delta_{p-2}= \frac{v^{-(p-3)}-v^{-(p-2)}v}{p},\\
& \delta_{p-3}= \frac{v^{-(p-4)}-v^{-(p-3)}v}{p},\\
& \vdots\\
& \delta_1 = \frac{1-v^{-1}v}{p},\\
\end{split}
\end{equation}
with $-p<\delta_i\leq 0$.
\begin{proof}
We start of  the relation in $\Z[G_p]$
\begin{displaymath}
P(\sigma)\times(\sigma-v)= v^{-(p-2)}\times \prod_{k=0}^{p-2} (\sigma-v^k)+p\times R(\sigma)\times(\sigma-v)=p\times Q(\sigma),
\end{displaymath}
with $Q(\sigma)\in\Z[G_p]$ because  $\prod_{k=0}^{p-2} (\sigma-v^k)=0$ in ${\bf F}_p[G_p]$ and so 
$\prod_{k=0}^{p-2} (\sigma-v^k)=p\times R_1(\sigma)$  in $\Z[G_p]$.
Then  we identify in $\Z[G_p]$ the  coefficients in the relation
\begin{displaymath}
\begin{split}
&(v^{-(p-2)}\sigma^{p-2}+v^{-(p-3)}\sigma^{p-3}+\dots+v^{-1}\sigma+1)\times(\sigma-v)=\\
&p\times (\delta_{p-2}\sigma^{p-2}+\delta_{p-3}\sigma^{p-3}+\dots+\delta_1\sigma+\delta_0),\\
\end{split}
\end{displaymath}
where $\sigma^{p-1}=1$.
\end{proof}
\end{lem}
\paragraph{Remark:}
\begin{enumerate}
\item
Observe that, with our notations,  $\delta_i\in \Z,\quad i=1,\dots,p-2$, but generally $\delta_i\not\equiv 0\modu p$. 
\item
We see also that $-p< \delta_i\leq 0$.
Observe also that $\delta_0=\frac{v^{-(p-2)}-v}{p}=0$.
\end{enumerate}
%
We shall show that the polynomial $Q(\sigma)$ caracterize in a straightforward  way the $p$-class group of $K_p$.
\begin{thm}\label{t512171}
Let $p$ be an odd prime. Let $v$ be a primitive root $\modu p$. 
For $k=1,\dots,r$ rank of the $p$-class group of $K_p$  then 
\begin{equation}\label{e512191}
Q(v^{2m_k+1})=\sum_{i=1}^{p-2} v^{(2m_k+1)\times i}\times(\frac{v^{-(i-1)}-v^{-i}\times v}{p})\equiv 0\modu p,
\end{equation}
(or an other formulation  $\prod_{k=1}^{r}(\sigma-v^{2m_k+1})$ divides $Q(\sigma)$ in ${\bf F}_p[G_p]$).
\begin{proof}$ $
\begin{enumerate}
\item
Let us fix $A$ one the singular numbers $A_k$ with $\pi^{2m+1}\  \|\ A-1$ equivalent to 
$\pi^{2m+1}\ \|\  (\frac{A}{\overline{A}} -1)$,
equivalent to 
\begin{displaymath}
\frac{A}{\overline{A}} =1+\lambda^{2m+1}\times a,\quad a\in K_p, \quad v_{\pi}(a)=0.
\end{displaymath}
Then raising to $p$-power we get 
$(\frac{A}{\overline{A}})^p=(1+\lambda^{2m+1}\times a)^p\equiv 1+p\lambda^{2m+1} a\modu\pi^{p-1+2m+2}$ and so 
$\pi^{p-1+2m+1}\ \|\  (\frac{A}{\overline{A}})^p -1.$
\item
From theorem \ref{l512164} we get 
\begin{displaymath}
(\frac{A}{\overline{A}})^{P(\sigma)\times(\sigma-v)} =(\frac{A}{\overline{A}})^{pQ(\sigma)}\equiv 1\modu\pi^{2(p-1)}.
\end{displaymath}
We have shown that 
\begin{displaymath}
(\frac{A}{\overline{A}})^p=1+\lambda^{p-1+2m+1}b,\quad b\in K_p,\quad v_\pi(b)=0,
\end{displaymath}
then 
\begin{equation}\label{e512281}
(1+\lambda^{p-1+2m+1}b)^{Q(\sigma)}\equiv 1\modu  \pi^{2(p-1)}.
\end{equation}
\item
But $1+\lambda^{p-1+2m+1}b\equiv 1+p\lambda^{2m+1} b_1\modu \pi^{p-1+2m+2}$ with $b_1\in\Z,\quad b_1\not\equiv 0\modu p$. 
There exists a natural integer $n$ not divisible by $p$ such that 
\begin{displaymath}
(1+p\lambda^{2m+1} b_1)^n\equiv 1+p\lambda^{2m+1}\modu\pi^{p-1+2m+2}.
\end{displaymath}
Therefore 
\begin{equation}\label{e512192}
(1+p \lambda^{2m+1}b_1)^{n Q(\sigma)}\equiv (1+p \lambda^{2m+1})^{Q(\sigma)}\equiv 1\modu  \pi^{p-1+2m+2}.
\end{equation}
\item
Show that the possibility of climbing up the  step  $\modu \pi^{p-1+2m+2}$  implies that $\sigma-v^{2m+1}$ divides $Q(\sigma)$ in ${\bf F}_p[G_p]$:
we have $(1+p\lambda^{2m+1} )^\sigma=1+p\sigma(\lambda^{2m+1}) = 1+p(\zeta^v-1)^{2m+1} =1+p((\lambda+1)^v-1)^{2m+1}
\equiv 1+pv^{2m+1}\lambda^{2m+1}\modu\pi^{p-1+2m+2}.$
In an other part $(1+p\lambda^{2m+1})^{v^{2m+1}}\equiv 1+p\lambda^{2m+1} v^{2m+1}\modu\pi^{p-1+2m+2}$.
Therefore  
\begin{equation}\label{e512282}
(1+p\lambda^{2m+1})^{\sigma-v^{2m+1}}\equiv 1\modu\pi^{p-1+2m+2}.
\end{equation}
\item
By euclidean division of $Q(\sigma)$ by $\sigma-v^{2m+1}$ in ${\bf F}_p[G_p]$, we get 
\begin{displaymath}
Q(\sigma) = (\sigma-v^{2m+1}) Q_1(\sigma)+R
\end{displaymath}
with $R\in{\bf F}_p.$
From congruence (\ref{e512192}) and (\ref{e512282})  it follows that  $(1+p\lambda^{2m+1})^R\equiv 1\modu\pi^{p-1+2m+2}$ and so 
that $1+pR\lambda^{2m+1}\equiv 1\modu\pi^{p-1+2m+2}$ and finally 
that $R=0$. 
Then in ${\bf F}_p$ we have $Q(\sigma)=(\sigma-v^{2m+1})\times Q_1(\sigma)$ and so  $Q(v^{2m+1})\equiv 0\modu p$, or explicitly 
\begin{displaymath}
\begin{split}
& Q(v^{2m+1})=v^{(2m+1)(p-2)}\times \frac{v^{-(p-3)}-v^{-(p-2)}v}{p}+v^{(2m+1)(p-3)}\times\frac{v^{-(p-4)}-v^{-(p-3)}v}{p}+\dots\\
&+v^{2m+1}\times \frac{1-v^{-1}v}{p}\equiv 0\modu p,\\
\end{split}
\end{displaymath}
which achieves the proof.
\end{enumerate}
\end{proof}
\end{thm}
\paragraph{Remarks:}
\begin{enumerate}
\item
Observe   that it is the $\pi$-adic theorem  \ref{l512164} connected to Kummer-Stickelberger which allows to obtain this result.
\item
It is easy to verify the consistency of relation (\ref{e512191}). 
The little following  MAPLE program allows to verify that this formula is consistent with the table of irregular primes and Bernoulli numbers in 
Washington, \cite{was} p. 410. We see that  $Q(v^{2m+1})\equiv 0\modu p$ when $B_{p-1-2m}\equiv 0\modu p$ where $B_{p-1-2m}$ is an even Bernoulli number.
\begin{verbatim}
restart;
> # p       : prime  in paper
> # delta_i :  in paper
> # Q_k     :   is Q in paper 
> # a_2     : is 2a in Washington table of irregular primes p. 410
> # v       : primitive root mod p  
> p:=3:
> p_max:=157:
> while p < p_max do :
>   p:=nextprime(p):
>   for v from 2 to p-2 do:
>   primitive:=1;
>   for i from 1 to p-2 do:
>   if (v^i mod p) =1 then primitive:=0:fi:
>   od:
>   if primitive = 1 then break:fi:
>   od:
> #***************
> for k from 3 to p-2 by 2 do:
> a2:=p-k;
>  Q_k:=0:
>  for i from 1 to p-2 do:
>   delta_i:=v^(-(i-1)) mod p;
>   delta_i:=delta_i-(v^(-i) mod p)*v:
>   #print(`delta_i=`,delta_i):
>   delta_i:=delta_i/p;
>   #print(`p=`,p,`v=`,v,`delta_i=`,delta_i):
>   Q_k:=Q_k+(v^(k*i) mod p)*delta_i :
>  od:
> Q_k:=Q_k mod p:
> if Q_k=0 then print(`p=`,p,`v=`,v,`a2=`,a2);fi;
> od:
> od:
> print (`****`);
                        p=, 37, v=, 2, a2=, 32
                        p=, 59, v=, 2, a2=, 44
                        p=, 67, v=, 2, a2=, 58
                        p=, 101, v=, 2, a2=, 68
                        p=, 103, v=, 5, a2=, 24
                        p=, 131, v=, 2, a2=, 22
                        p=, 149, v=, 2, a2=, 130
                        p=, 157, v=, 5, a2=, 110
                        p=, 157, v=, 5, a2=, 62
                                 ****

\end{verbatim}
\end{enumerate}
%
An  immediate consequence is an explicit  criterium for $p$ to be a regular prime:
\begin{cor}\label{512301}
Let $p$ be an odd prime. Let $v$ be a primitive root $\modu p$.
If the congruence  
\begin{equation}\label{e512301}
\sum_{i=1}^{p-2} X^{i-1}\times(\frac{v^{-(i-1)}-v^{-i}\times v}{p})\equiv 0\modu p
\end{equation}
has no solution $X$ in $\Z$ then the prime $p$ is regular.
\end{cor}
%
%

%
\end{document}